\documentstyle[12pt]{article}
\let\Bbb\bf
\def\be{\begin{equation}}
\def\ee{\end{equation}}
\def\bea{\begin{eqnarray*}}
\def\eea{\end{eqnarray*}}
\newtheorem{lemma}{Lemma}
\newtheorem{definition}{Definition}
\newtheorem{theorem}{Theorem}
\newtheorem{proposition}{Proposition}

\newenvironment{proof}{\medskip {\bf Proof.}}{\hfill \rule{.5em}{1em} \\}

\def\+{\oplus}

\def\*{^{\ast}}

\def\k{\kappa}
\def\s{\sigma}

\def\bp{{\Bbb P}}
\def\BP{{\Bbb CP}}

\def\br{{\Bbb R}}
\def\bc{{\Bbb C}}
\def\bcp{{\Bbb CP}}

\def\O{{\cal O}}
\begin{document}
\sloppy
\title{Twistors, K\"ahler Manifolds, and\\ Bimeromorphic Geometry\\ II}

\author{\parbox{2in}{\center Claude LeBrun\thanks{Supported
in part by  NSF grant DMS-9003263.}\\SUNY Stony
 Brook} ~~~
  {\em and}~~~
 \parbox{2in}{\center Yat-Sun Poon\thanks{Supported
in part by  NSF grant  DMS-8906806.}\\Rice University} }
\date{}
\maketitle

\begin{abstract}
Using examples \cite{13}
of compact complex 3-manifolds which arise as  twistor spaces,
we show that
the class of compact
complex manifolds  bimeromorphic to K\"ahler manifolds
is not stable under  small deformations of complex structure.
\end{abstract}
\vfill
 \pagebreak

A well-known theorem of Kodaira and Spencer \cite{11}\cite{15}
states that any small
deformation of the complex structure  of a compact  K\"ahler   manifold
 again  yields a complex  manifold
of K\"ahler type.
The question has been therefore been
 raised \cite{7} \cite{19}
as to whether a similar stability result holds for  compact  complex
manifolds
which are {\em bimeromorphically
equivalent\/}\footnote{Two connected compact
 complex m-manifolds $X$ and $Y$ are called {\em bimeromorphically
equivalent} if there exists a  complex m-manifold
 $V$, and   degree 1 holomorphic maps
  $V\to X$ and $V\to Y$.} to K\"ahler manifolds---
that is,  for manifolds of  Fujiki's class
$\cal C$ \cite{6}.
In this article, we will
analyze the twistor spaces obtained in the previous article \cite{13}
as small deformations of the Moishezon twistor
 spaces discovered in \cite{12}, and show that they
 are generically  {\em not} spaces
 of class $\cal C$, even though they are obtained as small
deformations of spaces which {\em are\/}. In short,
the bimeromorphic analogue of the Kodaira-Spencer stability theorem is
false.\footnote{A technically  different proof of this
 result, incorporating information exchanged in
 discussions and letters  with the first author
during the summer of 1990,
 was found  simultaneously by F. Campana \cite{3}, who has
chosen to publish his work separately.}

In an attempt to make this article as self-contained as possible,
 we begin with a brief
introduction to the subject, including a
quick review of the essential results of  the preceding article \cite{13}.

\bigskip

Our focus here will be on  the following
class of complex manifolds:

\begin{definition}
A {\em twistor space} will herein  mean a
 compact  complex 3-manifold $Z$  with the following properties:
\begin{itemize}
\item  There is a free anti-holomorphic involution $\s : Z\to Z$, $\s^2=$
identity,
called the {\em real structure} of $Z$;
\item There is a foliation of $Z$ by $\s$-invariant holomorphic curves
$\cong \bcp_1$, called the {\em real twistor lines}; and
\item Each real twistor line has normal bundle  holomorphically isomorphic to
$\O (1)\+\O (1)$, where $\O (1)$ is the degree-one line bundle on $\bcp_1$.
\end{itemize}
\end{definition}
\noindent
The space  $M$ of real twistor lines is thus a compact real-analytic
4-manifold,
 and we have  real-analytic submersion $\wp :Z\to M$ known as the
{\em twistor projection}. By a construction  discovered by Roger Penrose
\cite{16}, the complex structure of $Z$ induces
 a half-conformally-flat
conformal Riemannian conformal metric on $M$, and every such metric
conversely arises
in this way \cite{1};
 however, we will never explicitly need this  in the sequel.

We will only  concern ourselves here with the  class of
twistor spaces admitting hypersurfaces of the following type:

\begin{definition}
An {\em elementary divisor} $D$ on a twistor space $Z$ is a
complex hypersurface $D\subset Z$ whose homological intersection number
with a twistor line is +1, and such that
 $D\cap \s (D)\neq \emptyset$.
\end{definition}

An  elementary divisor is necessarily a smooth  hypersurface.
The existence of such a  divisor $D$  is a powerful
hypothesis indeed, for  it follows
(\cite{13}, Proposition 6)
that $D$ is an $n$-fold  blow-up
of $\bcp_2$, that $M$ is diffeomorphic to an
 $n$-fold connected sum
$\bcp_2\#\cdots\# \bcp_2$, and  the  map $\wp|_D:D\to M$
 contracts a projective line to a point, but is  elsewhere an
orientation-reversing  diffeomorphism.

In fact, these conclusions are quite sharp.

\begin{proposition} Let $X$ be any compact complex surface obtained from
$\bcp_2$ by blowing up distinct points. Then there exists a twistor space
$Z$ which contains an  elementary  divisor $D$  such that
$D\cong X$ as a  complex surface. Moreover, given a
smooth (respectively,  real-analytic) 1-parameter family
$X_t$ of surfaces obtained from
$\bcp_2$ by blowing up distinct  ordered
 points, there is a smooth (respectively,  real-analytic)
family $(Z_t, D_t)$ of twistor spaces
with  elementary divisors  such that $D_t\cong X_t$.\label{def}
\end{proposition}

\begin{proof} In \cite{12} it was shown that, given an arbitrary  blow-up $D$
of
$\bcp_2$ at $n$  collinear points, there is a twistor space $Z$ containing
a degree 1 divisor isomorphic to $D$.  In fact, such twistor spaces
$Z$ may be explicitly
constructed from   conic bundles  over $\bcp_1\times\bcp_1$ by
a process of blowing  subvarieties up and down, and thus
 may be taken to be {\em Moishezon} in this case. In the
accompanying article \cite{13}, the deformation theory of these twistor spaces
 was  studied, with the following conclusion.
Let $p_1:=(0,0)$ and $p_2:=(1,0)$ in $\bc^2$, and let
${\cal W}\subset [\bc^2]^{n-2}$ denote the set
$$\{ (p_3, \ldots , p_n)|p_j\in \bc^2, p_j\neq p_k, j,k=1,\ldots , n\};$$ let
${\cal L}\subset {\cal W}$ denote the subset
$p_3, \ldots , p_n\in (\bc \times \{ 0\})$ of collinear configurations.
 It was shown there
(\cite{13}, {\bf Theorem 3}) that there exists a (versal) family
$({\cal Z}, {\cal D})$ of twistor spaces with
elementary divisors over a ${\cal U}$ neighborhood of $[{\cal L}\times
(\br^+)^n]
\subset [{\cal W}\times (\br^+)^n]$ such that the
 divisor $D$ associated  with a  configuration of points
$$p_1,\ldots , p_n\in {\Bbb C}^2\subset\bcp_2$$ and arbitrary collection
of positive weights
$$m_1,\ldots , m_n\in \br^+$$
is isomorphic to
$\bcp_2$ blown up at $p_1,\ldots , p_n$.

Now suppose we are given an arbitrary compact complex surface
$X$  obtained from
$\bcp_2$ by blowing up $n$ distinct points $q_1, \ldots , q_n$. There is a
line $L\subset\bcp_2$ which misses $q_1, \ldots , q_n$, and now identify
$\bcp_2-L$ with $\bc^2$ in such a way that
$q_1=(0,0)$ and $q_2=(1,0)$. Assign all the points, say, weight 1.
By making a linear transformation, we may
also take the
points $q_1, \ldots , q_n$ to be as close as we like to the $z_1$-axis,
so that our configuration becomes a point of ${\cal U}$.
The corresponding fiber of our family $({\cal Z}, {\cal D})$ then comes
equipped with an  elementary divisor isomorphic to the given $X$.

On the other hand, suppose we are  instead given an arbitrary smooth
 family $X_t$ of surfaces obtained by
blowing up $n$ distinct, ordered points in $\bcp_2$, where $t$ ranges over
$\br$.
Let ${\cal X}\to \br$ denote the family with
fibers  $\{ X_t\}$. There is a bundle  ${\cal P}\to B$ of $\bcp_2$'s
from which ${\cal X}\to \br$ is obtained by blowing up $n$ sections
$q_1, \ldots , q_n$; let ${\cal P}^{\ast}\to \br$ denote the bundle of dual
planes,
in which the $q_1, \ldots , q_n$ define $n$ complex hypersurfaces.
The complement of these hypersurfaces in ${\cal P}^{\ast}$
has real codimension 2,  so, by transversality, a generic smooth
(respectively real-analytic)  section of ${\cal P}^{\ast}$ will
miss them, and we may therefore smoothly (respectively real-analytically)
choose a projective line $L_t$ in each fiber $P_t$ of ${\cal P}$
which misses the points ${q_1}_t, \ldots {q_n}_t$.
Using $q_1$ as the zero section, the complement of these chosen lines
becomes a vector bundle over $\br$ and so may be trivialized in such a manner
that $q_2\equiv (1,0)$. Our family of surfaces may therefore be thought of as
associated with a family of point configurations $(q_1, \ldots , q_n)_t$
in $\bc^2$,
where $q_1\equiv (0,0)$ and $q_2\equiv (1,0)$.
Again, let us assign each point a positive weight, say 1.
Now there is a positive real-analytic function $F(\zeta_3,\ldots , \zeta_n)$
such that a weighted configuration $((0,0,1), (0,1,1), (\zeta_3, \eta_3, 1),
\ldots , (\zeta_n, \eta_n, 1))$ is in ${\cal U}$ provided that
$$\sum | \eta_j|^2 < F(\zeta_3,\ldots , \zeta_n).$$
Setting  $(q_3, \ldots , q_n)_t = ((\zeta_3(t), \eta_3(t)),
\ldots , (\zeta_n(t), \eta_n(t)))$,  define
$$(p_1, \ldots , p_n)_t:=
\left[ \begin{array}{cc}1&0\\0&\sqrt{\frac{F(\zeta_3(t),
\ldots , \zeta_n(t))}{
(1+\sum|\zeta_j(t)|^2)}}
\end{array} \right](q_1, \ldots , q_n)_t~.$$ The family
  $((p_1,1) , \ldots , (p_n,1))_t$
of
weighted configurations  then takes values within the
parameter space ${\cal U}$ of the family $({\cal Z}, {\cal D})$.
 Pulling back $({\cal Z}, {\cal D})$ now yields the desired family of
twistor spaces with elementary divisors.
\end{proof}

It might be emphasized, incidentally, that the twistor space $Z$
is  by no means  determined by the intrinsic structure of a
 elementary divisor $D$.
Nonetheless, we will presently   see  that the intrinsic structure of such a
divisor {\em does} tell us a great deal about a twistor space, and is,
in particular, sufficient to determine its
{\sl algebraic dimension}.

Let us recall that the algebraic dimension $a(Z)$
 of a compact complex manifold $Z$ is by definition
the degree of transcendence its the field of meromorphic
 functions, considered as an extension of the field $\bc$ of constant
functions.
Equivalently, the algebraic dimension of $Z$ is precisely
the maximal possible dimension of
the image of $Z$ under a  meromorphic map to $\bcp_N$; in particular,
$a(Z)\leq \mbox{dim}_{\bc}(Z)$. When equality is
achieved in the latter inequality, $Z$ is said to be a
{\em Moishezon manifold\/} \cite{14}, and a suitable sequence of
blow-ups of $Z$ along complex submanifolds will then result in a
projective variety.

The following lemma of F. Campana will be of critical importance:

\begin{lemma} {\rm \cite{2}.} A twistor space $Z$ is
bimeromorphic to a K\"ahler
manifold iff it is Moishezon.\label{camp}
\end{lemma}
\begin{proof}
 Let $p$ and $q$  be
 distinct points of a real twistor line $L$ in a twistor space $Z$,
 and let $S_p$ (respectively, $S_q$) denote the space
of rational curves through $p$  (respectively, $q$)  which are
deformations of $L$.
Assume that  $Z$ is in the  class $\cal C$.
Because the  components of the
 Chow variety of $Z$ are therefore
compact, the correspondence space $$Z':=\{ (r, C_1, C_2)\in
Z\times S_p\times S_q~|~r\in C_1\cap C_2\}$$ is thus a
 compact complex space; by blowing up any singularities, we may assume that
$Z'$ is smooth.
But since  a real twistor
line has the same normal bundle as a projective line in $\bcp_3$,
a generic point of $Z$ is joined to  either $p$ or $q$ only by a discrete
set of curves of the fixed class. The correspondence space
$Z'$ is therefore generically a branched cover of $Z$, and is, in
 particular, a 3-fold.
On the other hand, we have a canonical map $$\phi:
Z'\to \bp (T_pZ)\times \bp (T_pZ)
\cong \bcp_2\times\bcp_2$$
obtained by taking the tangent spaces of
curves at their base-points $p$ or $q$.
Let $r$ be a point of $Z$ which is not on $L$, but close
enough to $L$ so that $r$ is joined to $p$ and $q$ by small deformations $C_1$
and $C_2$ of $L$, both of which  are  $\bcp_1$'s with normal bundle
$\O (1)\oplus\O (1)$. Then $(r,C_1, C_2)$ is a  point of
 $Z'$ at which the
 derivative of $\phi$ has maximal rank. Pulling back meromorphic functions
from  $\bcp_2\times\bcp_2$ to $Z'$ must therefore
yield 3 algebraically independent functions, and
$Z'$ is therefore a Moishezon space.
But since the projection $Z'\to Z$ is  surjective,
and since  the class of Moishezon manifolds is closed under holomorphic
surjections \cite{14}, it follows
that $Z$ is also a  Moishezon manifold.

The converse  is, of course, trivial. \end{proof}

On the other hand, the following lemma allows one to
determine the algebraic dimension of a twistor space:

\begin{lemma} {\rm \cite{17}.} Any meromorphic function on
a simply-connected twistor space $Z$ can be expressed as the ratio of two
holomorphic sections of a sufficiently large power $\k^{-m}$ of
the anti-canonical line bundle  $\k^{-1}:=\wedge^3TZ$. \label{poon}
\end{lemma}
\begin{proof}
We begin by observing that any
(compact)  twistor space satisfies $h^1(Z,\O )=b_1(Z)$.
This is a consequence of the {\em Ward correspondence\/} \cite{20}, which says
that
the set of holomorphic vector  bundles on $Z$ which are trivial
on real twistor lines  is in  1-1 correspondence with the
instantons on $M$; in particular,
every holomorphic line bundle on $Z$
with $c_1=0$ is obtained by pulling back a flat $\bc_{\ast}$-bundle from $M$
and equipping it with the obvious holomorphic structure.
With the exponential sequence
$$\cdots\to H^1(Z,\O )\to H^1(Z,\O_{\ast})\stackrel{c_1}{\to}
 H^2(Z,{\Bbb Z})\to \cdots$$
this implies that holomorphic
line bundles on a simply-connected twistor space are
classified by their Chern classes.

Since we have assumed that  $Z$ is simply connected, it follows that
$H^2(Z, {\Bbb Z})$ is free. On the other hand, the Leray-Hirsch
theorem tells us that $H^2 (Z, {\Bbb Q})={\Bbb Q}c_1(Z)\oplus H^2(M,{\Bbb Q})$.
The latter splitting of the cohomology is exactly  the
decomposition of $H^2 (Z, {\Bbb Q})$ into the $(\mp 1)$-eigenspaces
 of $\sigma^{\ast}$; a class will be called {\em real } if it is in the
$(-1)$-eigenspace, and a complex line-bundle will be called real if its
first Chern class is real.  There is thus a unique ``fundamental''
holomorphic line bundle  $\xi$ on $Z$ such that any real holomorphic
line bundle is a power of $\xi$ and such that the restriction of
$\xi$ to a twistor line is positive; in particular,  $\k=\xi^k$ for some $k$.
While we will not need to know this explicitly, it can in fact be shown
\cite{9} that $k=4$ if $M$ is spin, and $k=2$ otherwise.

Now suppose that we are given a meromorphic function $f$ on such a $Z$.
The function $f$ can {\em a priori} be  expressed in the
form $f=g/h$, where $g$ and $h$ are holomorphic sections of a line-bundle
$\eta\to Z$; for example, we could
take $\eta$ to be   the divisor line bundle of the polar locus of $f$.
The pull-back  $\sigma^{\ast}\overline{\eta }$ of the conjugate line-bundle
of $\eta$ is automatically holomorphic, and  $\sigma^{\ast}\overline{g}$
and $\sigma^{\ast}\overline{h}$ are holomorphic sections of this bundle.
The  holomorphic bundle $\eta\otimes\sigma^{\ast}
\overline{\eta}$ is now {\em real} and has sections, and so must be the form
$\xi^m$ for some positive integer $m$.
Thus
$$f=\frac{gh^{k-1}\sigma^{\ast}\overline{h}^k}{h^{k}\sigma^{\ast}\overline{h}^k},$$
expresses our meromorphic function as the quotient of two holomorphic sections
of $\k^m$.
\end{proof}

 We have already seen that there
are examples of Moishezon twistor spaces $Z$ containing an
elementary divisor $D$ isomorphic to $\bcp_2$ blown up at a collinear
configuration of points.
 We will now see that
the situation is dramatically different when the intrinsic structure of $D$ is
generic.

\begin{proposition}  Suppose that $Z$ is a twistor space with an  elementary
divisor isomorphic to the blow-up of
$\bcp_2$   at
$n$ generic points, $n\geq 7$. Then  $Z$ has no
non-constant meromorphic functions, and so has algebraic dimension $0$.
The set of configurations $(p_1 ,\ldots ,p_n )$ of points in $\bcp_2$
which are generic in this sense
is the complement of a countable union of proper algebraic sub-varieties
of $(\bcp_2)^n$, and in particular has full measure.
\label{big}
\end{proposition}
\begin{proof}
Let us begin by considering  the case of a
 configuration of $n$ points in ${\Bbb C}^2$ containing a 6-point configuration
of the following type:

\setlength{\unitlength}{3ex}
\begin{center}\begin{picture}(18,5)(0,0)
\put(0,0){\line(3,1){12}}
\put(18,0){\line(-3,1){12}}
\put(0,1){\line(1,0){18}}
\put(3,1){\circle*{.2}}
\put(9,1){\circle*{.2}}
\put(9,3){\circle*{.2}}
\put(6,2){\circle*{.2}}
\put(12,2){\circle*{.2}}
\put(15,1){\circle*{.2}}
\end{picture}
\nopagebreak

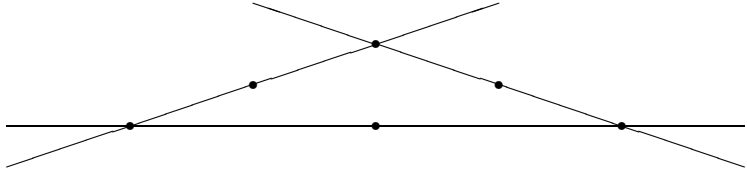
{\bf  Figure 1.} \end{center}
We assume that the other points of the configuration are
 not on any of three projective lines
of the figure.
The proper transforms of  these three  lines
are then (-2)-curves $E_j$, $j=1,2,3$. The anti-canonical bundle
$\kappa_D^{-1}$ of the surface $D$ thus
 satisfies $\kappa_D^{-1}|_{E_j}\cong {\cal O}$.
On the other hand, since the
half-anti-canonical bundle of $Z$ is given by
$\kappa^{-1/2}= [D] \otimes [\overline{D}]$, we have
$$\kappa^{-1/2}|_{D}=\nu\otimes [L_{\infty}]~,$$
where $\nu$ denotes the normal bundle of $D\subset Z$ and
$L_{\infty}\subset D$ is the projective line $D\cap \overline{D}$.
Yet the adjunction formula yields
$$\kappa^{-1}|_{D}=\nu \otimes \kappa_D^{-1}~,$$
so that
$$\nu^{2}\otimes [L_{\infty}]^2= \nu \otimes \kappa_D^{-1}$$
implying that
$\nu = \kappa_D^{-1}\otimes [L_{\infty}]^{-2}$
and hence
\be \kappa^{-1/2}|_{D}= \kappa_D^{-1}\otimes [L_{\infty}]^{-1}~.\label{snuk}\ee
It follows that
$$\kappa^{-1/2}|_{E_j}\cong {\cal O}(-1)~.$$
On the other hand, the normal bundle $N_j$ of $E_j\subset D$ is isomorphic to
${\cal O}(-2)\to {\Bbb CP}_1$. Since
\bea
\Gamma (E_j, {\cal O} ((\kappa^{-m/2}|_{E_j})\otimes N_j^{-k}))&=&
\Gamma ({\Bbb CP}_1, {\cal O} (-m+2k))\\&=& 0 ~~\makebox{if $k<
\frac{m}{2}$,}\eea
it follows that any section of $\kappa^{-m/2}|_{D}$ vanishes along
$E_j$ to order $[\frac{m-1}{2}]$. But through the generic
 point of $D$ we can find
a projective line  in  $D$ passing through  a
 blown-up point not on  the diagram, avoiding all other blown-up points,
and meeting the $E_j$ in three distinct points.
 Letting $L$ denote the proper transform of such a line, one has
\bea
\kappa^{-1/2}|_L&=& (\kappa_D^{-1}\otimes [L_{\infty}]^{-1})|_L\\&\cong&
{\cal O}(2)\otimes {\cal O}(-1)\\&\cong&{\cal O}(1)~,\eea
so that $\kappa^{-m/2}|_L\cong {\cal O}(m)$. Yet any holomorphic section of
of $\kappa^{-m/2}|_{D}$ must have 3 zeroes on $L$ of multiplicity
$[\frac{m-1}{2}]$ at $L\cap E_j$.  Since $3[\frac{m-1}{2}]> m$ for $m>6$,
we conclude that such a section must vanish identically on $L$ provided $m$
is sufficiently large. Hence $\Gamma (D, {\cal O} (\kappa^{-m/2}))=0$
for $m$ sufficiently large, and hence, by taking tensor powers
of sections, for all $m>0$. Similarly,
 $\Gamma (\overline{D}, {\cal O} (\kappa^{-m/2}))=0$ for all $m>0$. From the
exact sequences
\be 0\to  {\cal O}_{Z} (\kappa^{-(m-1)/2})\to {\cal O}_{Z}  (\kappa^{-m/2})\to
 {\cal O}_{D\cup
\overline{D}}  (\kappa^{-(m-1)/2})\to 0,\label{snork}\ee
 we conclude by induction that
$$\Gamma ({Z}, {\cal O}(\kappa^{-m/2}))=\bc$$
for all $m> 0$. By Lemma \ref{poon},
 any meromorphic function on $Z$ must therefore be constant.

We now examine the case of $D$ obtained from $\bcp_2$ by blowing up
$n>6$ generically located points. For each n-tuple of points  $(p_1 ,\ldots ,
p_n )_u$ in $\bc^2=\bcp_2-L_{\infty}$, let $D_u$ denote the
corresponding blow-up of $\bcp_2$, and consider  the behavior
of $h^0(D_u , {\cal O} (\kappa_{D_u}^{-m}\otimes [L_{\infty}]^{-m}))$.
By the semi-continuity principle \cite{8}
and the above calculation, this vanishes, for $m$ fixed,  on a
non-empty Zariski-open subset of configurations. The set of $n$-point
configurations for which
$h^0(D_u , {\cal O} (\kappa_{D_u}^{-m}\otimes [L_{\infty}]^{-m}))\neq 0$
for some $m$
is therefore a countable union of subvarieties, and so
has measure $0$. Using the exact sequence \ref{snork} and the isomorphism
\ref{snuk}, we conclude that
$\Gamma ({Z}, {\cal O}(\kappa^{-m/2}))=\bc ~\forall m\neq 0$
provided that $Z$ contains an elementary divisor $D$ obtained from
$\bcp_2$ by blowing up $n>6$ generic points.
 Again applying Lemma \ref{poon}, we conclude that, for $n\geq 7$,
 any meromorphic function on a twistor space
$Z$ containing  a generic elementary divisor must therefore be constant.
\end{proof}

 Our main result now follows  immediately:

\begin{theorem} The class $\cal C$, consisting of
 of compact complex manifolds which are
bimeromorphic to K\"ahler manifolds,
is  not stable under small deformations.
\end{theorem}
\begin{proof} By Propositions \ref{def} and \ref{big}, there exist
1-parameter families of twistor spaces  $Z_t$ for which almost every
$Z_t$ has algebraic dimension $0$, whereas $Z_0$ is Moishezon;
in fact it suffices to take $Z_0$ to be one of the explicit examples of
\cite{12},
with $D_0$ corresponding to a collinear configuration of
$n\geq 7$ points, arrange for the
curve of configurations $(p_1, \ldots ,
p_n)_t$  to  be real-analytic and
contain at least one generic configuration. (Actually,
one can do better:
by taking the elementary divisors $D_t$ to all correspond to
configurations containing projective copies of Figure 1
when $t\neq 0$, one can even arrange  for
 $Z_0$ to be the {\em only} Moishezon space in the
family.)
By   Lemma \ref{camp}, the non-Moishezon twistor spaces of the
family $Z_t$ are not of class $\cal C$, despite the fact that
 they are arbitrarily small deformations of the Moishezon space $Z_0$.
\end{proof}

\noindent {\bf Remarks.}
\begin{itemize}
\item
In order to keep this article as short and clear as possible,
we have only considered the case of $n\geq 7$, and only presented the
extreme cases of $a(Z)=3$ and $a(Z)=0$. In fact, in can be shown that
generically  $a(Z)<3$    as soon as $n\geq 4$. One can also find
simple non-collinear  configurations for which $n=1, 2$ as soon as
$n\geq 5$. Finally, one can show that the existence
of an elementary divisor corresponding to a collinear configuration
{\em forces} $Z$ to be one of the examples of \cite{12}, and, in particular,
 Moishezon.
For details, see \cite{18}.
\item The existence of self-dual metrics on
arbitrary connected sums $\BP_2\# \cdots\#
\BP_2$  was first proved
abstractly  by Donaldson and
Friedman \cite{4} and, using completely different methods,
by Floer \cite{5}. Unlike the methods used  here,  these
methods do not show that the twistor space of some of these
metrics are Moishezon.
It was nonetheless the Donaldson-Friedman construction  which  originally gave
the authors reason to believe that the generic deformation of
the  explicit twistor spaces of \cite{12}
should  not be of Fujiki-class $\cal C$.
For providing  this source of inspiration, as well as
for their friendly advice and encouragement, the authors would therefore like
to gratefully thank Robert Friedman and Simon Donaldson.
\end{itemize}

\end{document}